\documentclass{ifacconf}

\usepackage{graphicx}      
\usepackage{color}
\definecolor{hl}{rgb}{1, 0.1, 0.1}
\usepackage{tabulary}
\usepackage{multirow}
\usepackage{tabularx} 
\newcolumntype{C}{>{\centering\arraybackslash}X}
\newcommand{\fref}[1]{(\ref{eq:#1})}
\usepackage{amsfonts}
\usepackage{amsmath}
\usepackage{float}

\usepackage{natbib}        
\begin{document}
\begin{frontmatter}

\title{Market-Driven Energy Storage Planning for Microgrids with Renewable Energy Systems Using Stochastic Programming
}


\author{Abdulelah H Habib}, 
\author{Vahid R. Disfani}, 
\author{Jan Kleissl},
\author{Raymond de Callafon} 


\address{Mechanical and Aerospace Engineering\\ University of California San Diego,
La Jolla, CA 92093\\
\{ahhabib, disfani, jkleissl, callafon\}@ucsd.edu\\}

\def\zack{\textcolor{green}}
\def\vahid{\textcolor{cyan}}
\def\habib{\textcolor{red}}

\begin{abstract}                
Battery Energy Storage Systems (BESS) can mitigate effects of intermittent energy production from renewable energy sources and play a critical role in peak shaving and demand charge management. To optimally size the BESS from an economic perspective, the trade-off between BESS investment costs, lifetime, and revenue from utility bill savings along with microgrid ancillary services must be taken into account. The optimal size of a BESS is solved via a stochastic optimization problem considering wholesale market pricing. A stochastic model is used to schedule arbitrage services for energy storage based on the forecasted energy market pricing while accounting for BESS cost trends, the variability of renewable energy resources, and demand prediction. The uniqueness of the approach proposed in this paper lies in the convex optimization programming framework that computes a globally optimal solution to the financial trade-off solution. The approach is illustrated by application to various realistic case studies based on pricing and demand data from the California Independent System Operator (CAISO). The case study results give insight in optimal BESS sizing from a cost perspective, based on both yearly scheduling and daily BESS operation.
\end{abstract}


\end{frontmatter}

\section{Introduction}
\def\jan{\textcolor{blue}}
\def\raymond{\textcolor{red}}
\def\zack{\textcolor{green}}
\def\habib{\textcolor{cyan}}
\def\delete{\textcolor{white}}

 
The need for a Battery Energy Storage System (BESS) to serve as a buffer for electric energy is palpable for microgrid systems that have a large penetration of intermittent renewable energy sources. A BESS may be economical for both islanded microgrids and a for grid-connected system, as a BESS increases reliability during outages and provides revenue or grid services such as peak shaving, voltage regulation, and arbitrage power trading during normal operation \citep{Lasseter,Donadee1,storagereview} .

Applications of a BESS can be found in various settings to assist with renewable power integration. It has been applied to the problem of harmonic distortion, generally known as voltage regulation, which may occur in standalone operation (islanding) of a microgrid \citep{batteruvoltage,sandiavoltage}. 
Specifically, a BESS reduces (or smoothens out) variability in Photo Voltaic (PV) and wind energy production \citep{teleke2010rule,zheng2015optimal} by different control strategies such as a rule-based control and a model predictive control (MPC). 
A BESS in conjunction with PV and demand forecasting can help shift renewable generation to times of higher power demand or lower electricity price via an MPC technique \citep{sevilla2015advanced}. 
A mathematical model for a large BESS system was performed in \citep{statemodelbattery} as a reduced four state space equations to model the relation between the bulk power grid and a BESS.

The benefits of BESS in coping with variable renewable energy production are evident, but the costs associated with financing and installing BESS are often prohibitive. A good example is a residential setting \citep{batteryresidential}, where a BESS may not produce sufficient revenue from energy arbitrage to achieve investment payback with the current BESS prices without government incentives. At the same time, BESS costs are anticipated to drop in the near future and investment banks are expecting the payback time for unsubsidized investment in electric vehicles (EV) combined with rooftop solar and BESS \citep{houchoisglobal} to reduce to around six to eight years. 
Also, the economies of scale due to the adoption of EV and rapid improvement of battery technologies will likely reduce BESS prices. The projected drop in pricing of BESS is also expected to lead to an additional three-year reduction in payback time to three years by 2030 \citep{nykvist2015rapidly,bankstudy}. 

Optimal BESS sizing from an economical perspective must find the optimal trade-off between critical design parameters that include BESS sizing, BESS life expectancy due to battery degradation and total revenue from utility bill savings due to energy arbitrage. Holistic BESS scheduling models that aim to capture all cost aspects were developed in \citep{nguyen2012new} to maximize the overall profit of an existing wind-storage system.  
Economic models were used in \citep{ornelas2014optimized} to predict the market price to optimize the operation of existing energy resources in a microgrid, but no future investments were considered. Operational stochastic control and optimization in \citep{zachar2016economic} were designed as an MPC to ensure sufficient energy as an economic dispatch problem.


Motivated by the need to find the optimal BESS investment as a function of time considering capital and O\&M costs, as well as operational revenues, this paper proposes a stochastic optimization approach that leverages mixed integer and real (convex) optimization to formulate financially optimal BESS sizing solutions. The stochastic optimization is used to address the variability in prediction and forecasting of energy and BESS pricing to determine when is the optimal time to invest in a BESS. The convex optimization is used to compute globally optimal solutions for BESS sizing parameters, given the operational model and the price variability in the day-ahead market. 

The paper is outlined as follows. First, the problem formulation and the system topology for financial optimization are summarized in Section~\ref{sec:systop}. The mathematical framework is summarized in Section~\ref{sec:schprog}, explaining the optimization techniques, objective functions and the constraints. Different operating scenarios are discussed and compared in Section~\ref{sec:casestudy} to cover cases of extreme high/low power variability in solar, wind and demand patterns. In Section~\ref{sec:results}, different BESS installation cases and optimal BESS sizing for a case study of a real microgrid are presented.

\section{System Topology {and Pricing}}\label{sec:systop}
\subsection{Microgrid and Market Structures}
\begin{figure}[ht]
\includegraphics[width=\columnwidth]{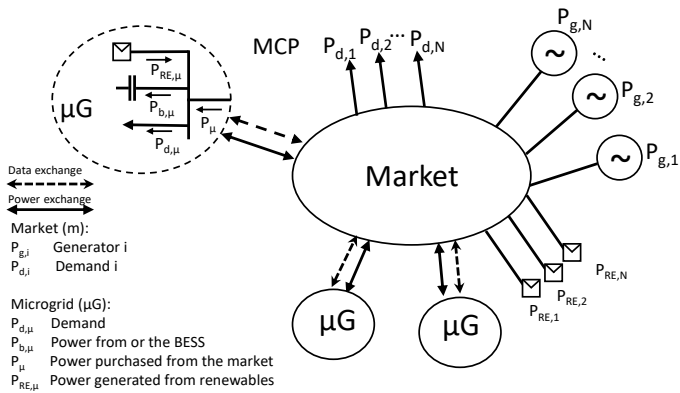}
\centering
\caption{Power market system architecture.}
\label{figure:system}
\end {figure}

Fig. \ref{figure:system} illustrates the structure of power market and microgrids used in this paper. The microgrid is modelled as a subset of the market $\mu G \subset \mathcal{N}$, and demand, renewable generation, and BESS power in both market ($m$) and microgrid ($\mu$) are denoted by $P_{d}$, $P_{RE}$, and $P_b$, as illustrated in Fig. \ref{figure:system}.

The net demand $P_{net}$, which is the actual market demand {(including all microgrids' demand)} 
minus the total renewable power available in the market, is computed via 
\begin{equation}\label{eq:price2}
P_{net}=\Sigma_{i\in \mathcal{N}}\left(P_{d_i} -P_{RE_i}\right)
=\Sigma_{i\in \mathcal{N}}P_{g_i}
\end{equation}
where $P_{g_i}$ indicates power provided by generator $i$ and
\begin{equation}\label{eq:price3}
P_\mu =\Sigma_{j\in \mu G}\left(P_{d_j,\mu}+P_{b_j,\mu} -P_{RE_j,\mu}\right)
\end{equation}

\subsection{Market Clearing Price Modeling}
Assuming that the microgrids will pay the hourly market clearing price (MCP) in the future instead of predefined constant or time-of-use (TOU) rates, a price model is required to anticipate the MCP at different times for optimal operation of microgrids. 

Typically, Independent System Operators (ISO) aggregate the bids received from generators and cross it with the net demand hourly profile of the market to define hourly MCP. It is assumed that the MCP is solely a function of the net demand in that {$\lambda_P$} is {linearly} correlated with the market net demand {$P_{net}$ via}
\begin{equation}
\lambda_P = \alpha~P_{net} +\beta\label{eq:price1}.
\end{equation}
This pricing modeling has been validated in the literature \citep{marketPrice1,marketprice2}. It is assumed for simplicity that the parameters of the MCP model ($\alpha, \beta$) remain constant throughout the 
15 year modeling horizon. However, the optimization could consider more detailed and dynamic models where the pricing model parameters vary as generators are added or removed.

To model the effects of different generators' bidding strategies and maintenance schedules on different days of the week (weekdays and weekend) and different seasons (summer and non-summer) on MCP, four distinct MCP models are fit from historical CAISO demand and pricing data.

\subsection{Microgrid and Power Market Growth}
For realistic financial predictions and optimal sizing of the BESS, the financial model considers the annual growth of both the market and the microgrid. The growth of the market and the microgrid takes into account all components, i.e. demand, solar and wind.

For simplicity, we assume a fixed annual solar growth (ASG) defined by
\[ASG=\frac{S_{t+1y}-S_t}{S_t}\times 100\%\]
where $S_t$ represents the vector of hourly solar profiles of the current year. 
Hence, with a fixed ASG, the net solar power $S_{t+1}$ contribution is predicted to grow exponentially as
\[
S_{t+1} = \left( \frac{ASG}{100} +1 \right ) S_t
\]
with $ASG>0$. Similarly, we assume a fixed annual wind growth (AWG) as
\[AWG=\frac{W_{t+1y}-W_t}{W_t}\times 100\%.\]
The CAISO historical demand data shows different rates of increase at different hours of the day, specifically power demand at the peak hour has grown faster than at off-peak hours. to account for this effect, we define an annual demand growth profile (ADGP) that varies by hour of the day as 
\[
ADGP=(D_{t+1}-D_t) \oslash {D_t}\times 100\%,\] where $D_t$ is the vector of hourly power demand at year $t$ and $\oslash$ denotes element-wise division. 

\section{Stochastic BESS Optimization}
\label{sec:schprog}

\subsection{Objective Function}
A stochastic programming model is developed to capture all possible scenarios of solar and wind generations, demand variations, 
and supply bidding strategies
in the BESS sizing and operation problem.

The objective function is defined to minimize the expected energy cost, i.e. the sum of energy purchase costs and the BESS investment cost over a time period $T$ while considering annual growth in solar, wind and demand. Mathematically, for an optimization horizon $T$ and scenario set of $\Omega$, the objective function is defined as
\begin{equation}
\\ \min_{P_b, {\bar{P}}_b}{} \Sigma_{t\in T}\left\{J_{b,t}+\Sigma_{i\in\Omega} Pr_i J_{i,t} -v_{b,T}\right\}\cdot v_t~, 
\label{eq:min}
\end{equation}
where $Pr_i$ refers to the probability of a scenario $i$ and 
\begin{itemize}
\item[i] $J_{b,t}$ refers to the investment cost of BESS installed in year $t$. 
\item [ii] $J_{i,t}$ is the total cost of energy purchased from the market in year $t$ for the scenario $i$, and is given by
\[J_{i,t}=\Lambda_{i,t}^T\rm E_{i,t}.\]
$\rm E_{i,t}$ and
$\Lambda_{i,t}$ are respectively the hourly vectors of amount of energy purchased from the energy markets and MCP in year $t$ and scenario $i$.
\item[iii] $v_{b,T}$ accounts for the remaining value of the unexpired BESS at the end of the simulation period.
\item[iv] As is common in economic models $v_t=(1-\gamma)^t$ discounts the monetary value in future years using an annual interest rate of $\gamma$.
\end{itemize}



\subsection{Constraints}

\subsubsection{Resource Adequacy:}
At any time step, the microgrid control center (MGCC) must ensure that there is adequate power to supply demand. 
Therefore, any mismatch between the power demand and the summation of solar power output, and BESS discharging power must be purchased from the market to keep the power balance at any time step $t$ and scenario $i$. This constraint is summarized as
\begin{equation}
\begin{array}{c}
P_{RE_{i,t}}+P_{b_{i,t}}+P_{{i,t}}=L_{i,t}
\end{array}
\label{eq:resource}
\end{equation}
where $P_{i,t}$ is the power purchased from the market.

\subsubsection {Battery Constraints:}


First, the battery charging/ discharging power must be between the limits, \textit{i.e.}
\[
-P_{b_{max}} \le {P_{b_{i,t}}}\le P_{b_{max}} 
\] 

To avoid damages due to a deep (dis)charge cycle of the battery, the stored energy in the battery is constrained by its maximum and minimum SOC limits ($\rho_{min},\rho_{max}$) as
\[ \rho_{min} { E}_{b_{max}} \le { E}_{b_{i,t}} \le \rho_{max} { E}_{b_{max}}\]
where $\rho_{min}$ and $\rho_{max}$ are typically around 10\% and 95\%. The energy stored in the battery is denoted by ${E}_{b_{i,t}}$ and
calculated via 
\[E_{b_{i,t}} = \Sigma_{h=1}^{t}P_{b_{i,h}} \Delta t +E_{b_{i,0}}\] with  ${{\rm E}}_{b_{i,0}}$ as the initial BESS energy and $\Delta t$ time difference between two consecutive time steps.

It is also desired to keep the final SOC of the BESS equal to  its initial value at the end of each day. This constraint is needed to avoid transferring energy between days and included via
\[{E_{b_{i,t_1}}}={E_{b_{i,0}}}\]for any ${t_1\in 
\{24k \text{ hours},~ k \in  \mathbb{N} \}}$.
Finally, the ratio between the nominal power rating and energy rating of the BESS implemented by 
\[
2\times { E}_{b_{max_i}} = P_{b_{max_i}}
\]
as the last battery constraint. Obviously, more advanced battery constrainst that take into account parasitic loss and efficiency parameters could be used to provide even more realistic battery constraints.

\subsubsection {Power Congestion Constraint}
A power congestion constraint limits the power purchased from the market due to the physical limit of the microgrid at the point of common coupling (PCC) or upstream power lines. Power congestion constraints may, for example, limit the BESS' ability to charge at $P_{b_{max_i}}$ during (or within) the cheapest price. By including a power congestion constraint
\begin{equation} \label{eq:Plimit}
-P_L \le P  \le P_L
\end{equation}
the BESS will be charged over a longer time frame to accommodate the congestion limit $P_{L}$.

\subsection{Scenarios}
\label{sec:Scenarios}

The most accurate results would be obtained by simulating a typical meteorological and climatologically representative timeseries over a year (or longer), but this is computationally intensive. Instead, we consider  year's (2015) worth of data, downsampled to a set of typical days and these days are repeated each year. 

\begin{figure}[ht]
\includegraphics[width=.85\columnwidth]{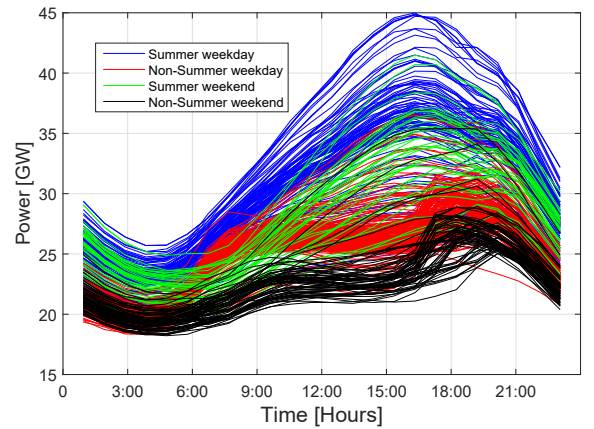}
\centering
\caption{Daily demand profiles for the CAISO market for one year (2015). Colors denote different clusters classified by (week)days and (non)summer season.}
\label{figure:loadprofile}
\end {figure}

\begin{figure}[ht]
\includegraphics[width=.85\columnwidth]{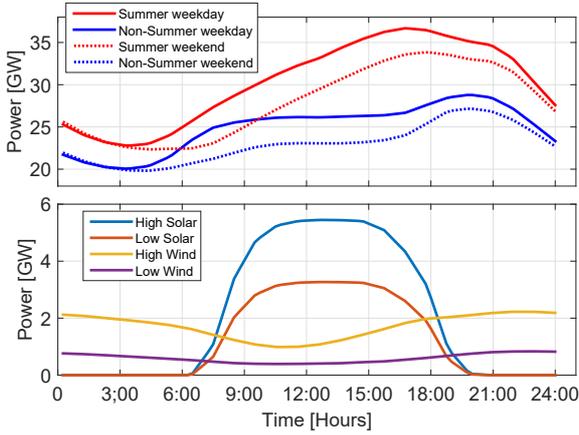}
\centering
\caption{Top: average of demand profiles clustered in 4 categories. Bottom: the 4 categories of high/low solar/wind and used as $\lambda_P$ in \fref{price1} for pricing modeling.}
\label{figure:Scenarios}
\end {figure}

Subsequently, a manual clustering is applied to classify demand profiles in representative patterns. The clustering assembles the data into four main groups that resemble a summer/non-summer and weekday/weekend separation. Fig. \ref{figure:loadprofile} illustrates all demand profiles clustered in those four groups, each identified by a distinct color. The clustering can then be used to formulate an average for each group as depcited in the top plot of Fig.~\ref{figure:Scenarios}.

Although usually intermittent in nature, solar and wind power profile over the complete state of California tends to be smooth and easily seperable into a small number of distinctive patterns. Only a binary classification of \textit{high} and \textit{low} is applied in this paper to cover the state-wide range of patterns in solar and wind power generation. Clusters were obtained by the popular clustering method k-means \citep{kmeans}, which aims to partition time series data into two clusters. Clear days are presented as the \textit{high} solar case, whereas overcast days are denoted by \textit{low} solar case. The bottom plot of Fig.~\ref{figure:Scenarios} illustrates the 4 different scenarios for high/low solar and wind.

For unambiguous notation we use \textit{H} and \textit{L} as abbreviated notations for the \textit{high} and \textit{low} and summer weekday and summer weekend and non-summer weekday and summer weekend as \{\textit{SWD,SED,NSWD,NSED}\} for power/demand conditions, while \textit{S}, \textit{W}, and \textit{D} are used for \textit{solar}, \textit{wind} and \textit{demand}. The binary classification with the demand scenarios allows $4\cdot 2\cdot2=~16$ scenarios. A summary of the 16 scenarios based on the binary classification of \textit{high} and \textit{low} solar, wind and demand
\begin{equation}
\begin{aligned}
\Omega =&\left\{HS,LS\right\}\times \left\{ HW, LW\right\}\\
&\times \left\{ SWD,SED,NSWD,NSED \right\}
\end{aligned}
\end{equation}

is given in Table~\ref{table:scenarios}. Clearly, more granular clustering, e.g. by adding seasonal effects on variables such as (solar, wind), will increase the accuracy of the optimization, but also the computational cost. 
\begin{table}[ht]
\scriptsize
\centering
\bgroup
\def\arraystretch{1.25}
\setlength\tabcolsep{-.2cm}
\caption{Scenarios $i$ for CAISO demand and renewable generation. The probability of each scenario $Pr_i$ is given in the last row.
}
\label{table:scenarios}
\begin{tabularx}{\columnwidth}{|C|C|C|C|C|C|C|C|C|C|C|C|C|C|C|C|C|}
\hline
 \parbox[t]{2mm}{\multirow{3}{*}{\rotatebox[origin=c]{90}{Demand}}}
& \multicolumn{8}{c|}{{\normalsize Summer}}                                                                                                                                           & \multicolumn{8}{c|}{{\normalsize Non-Summer}}                                                                                                                                              \\
                        & \multicolumn{4}{c}{\multirow{2}{*}{\normalsize weekday}}                                     & \multicolumn{4}{c|}{\multirow{2}{*}{\normalsize weekend}}                                     & \multicolumn{4}{c}{\multirow{2}{*}{\normalsize weekday}}                                        & \multicolumn{4}{c|}{\multirow{2}{*}{\normalsize weekend}}                                         \\
                        & \multicolumn{4}{c|}{}                                                             & \multicolumn{4}{c|}{}                                                             & \multicolumn{4}{c|}{}                                                                & \multicolumn{4}{c|}{}                                                                 \\ \hline
                         \parbox[t]{2mm}{\multirow{2}{*}{\rotatebox[origin=c]{90}{Solar}}}
                        & \multicolumn{2}{c|}{\multirow{2}{*}{H}} & \multicolumn{2}{c|}{\multirow{2}{*}{L}} & \multicolumn{2}{c|}{\multirow{2}{*}{H}} & \multicolumn{2}{c|}{\multirow{2}{*}{L}} & \multicolumn{2}{c|}{\multirow{2}{*}{H}}  & \multicolumn{2}{c|}{\multirow{2}{*}{L}}   & \multicolumn{2}{c|}{\multirow{2}{*}{H}}   & \multicolumn{2}{c|}{\multirow{2}{*}{L}}   \\
                        & \multicolumn{2}{c|}{}                   & \multicolumn{2}{c|}{}                   & \multicolumn{2}{c|}{}                   & \multicolumn{2}{c|}{}                   & \multicolumn{2}{c|}{}                    & \multicolumn{2}{c|}{}                     & \multicolumn{2}{c|}{}                     & \multicolumn{2}{c|}{}                     \\ \hline
                         \parbox[t]{2mm}{\multirow{2}{*}{\rotatebox[origin=c]{90}{Wind}}}& \multirow{2}{*}{H} & \multirow{2}{*}{L} & \multirow{2}{*}{H} & \multirow{2}{*}{L} & \multirow{2}{*}{H} & \multirow{2}{*}{L} & \multirow{2}{*}{H} & \multirow{2}{*}{L} & \multirow{2}{*}{H} & \multirow{2}{*}{L}  & \multirow{2}{*}{H}  & \multirow{2}{*}{L}  & \multirow{2}{*}{H}  & \multirow{2}{*}{L}  & \multirow{2}{*}{H}  & \multirow{2}{*}{L}  \\
                        &                    &                    &                    &                    &                    &                    &                    &                    &                    &                     &                     &                     &                     &                     &                     &                     \\ \hline
                         \parbox[t]{2mm}{\multirow{2}{*}{\rotatebox[origin=c]{90}{Scen.}}}                         & \multirow{2}{*}{1} & \multirow{2}{*}{2} & \multirow{2}{*}{3} & \multirow{2}{*}{4} & \multirow{2}{*}{5} & \multirow{2}{*}{6} & \multirow{2}{*}{7} & \multirow{2}{*}{8} & \multirow{2}{*}{9} & \multirow{2}{*}{10} & \multirow{2}{*}{11} & \multirow{2}{*}{12} & \multirow{2}{*}{13} & \multirow{2}{*}{14} & \multirow{2}{*}{15} & \multirow{2}{*}{16} \\                        &                    &                    &                    &                    &                    &                    &                    &                    &                    &                     &                     &                     &                     &                     &                     &                     \\ \hline
                                                \parbox[t]{2mm}{\multirow{3}{*}{\rotatebox[origin=c]{90}{~~~Pr}}}& \multirow{2}{*}{3.2} & \multirow{2}{*}{7.4} & \multirow{2}{*}{4.7} & \multirow{2}{*}{11.0} & \multirow{2}{*}{1.2} & \multirow{2}{*}{2.8} & \multirow{2}{*}{1.8} & \multirow{2}{*}{4.1} & \multirow{2}{*}{5.4} & \multirow{2}{*}{12.7} & \multirow{2}{*}{8.1} & \multirow{2}{*}{19.0} & \multirow{2}{*}{2.2} & \multirow{2}{*}{5.2} & \multirow{2}{*}{3.4} & \multirow{2}{*}{7.8} \\                        &                    &                    &                    &                    &                    &                    &                    &                    &                    &                     &                     &                     &                     &                     &                     &                     \\ \hline
             
\end{tabularx}
\egroup
\end{table}

Each scenario is assigned a probability consistent with climatological data in a certain location. The probability of each scenario is defined by multiplying the corresponding individual probabilities and is used as a weighting of that scenario in the optimization of \fref{min}, e.g.
\begin{equation}\label{eq:Pr}
Pr(\left\{HS, LW, SWD\right\})=Pr(HS)\cdot Pr(LW)\cdot Pr(SWD).
\end{equation}
Individual probabilities of solar, wind and demand are mutually independent and the probability of all possible scenarios sums to 100\%.

\section{Case Studies and Simulations }
\label{sec:casestudy}
\subsection{Market Clearing Price Models}
\begin{figure}[ht]
\includegraphics[width=.9\columnwidth]{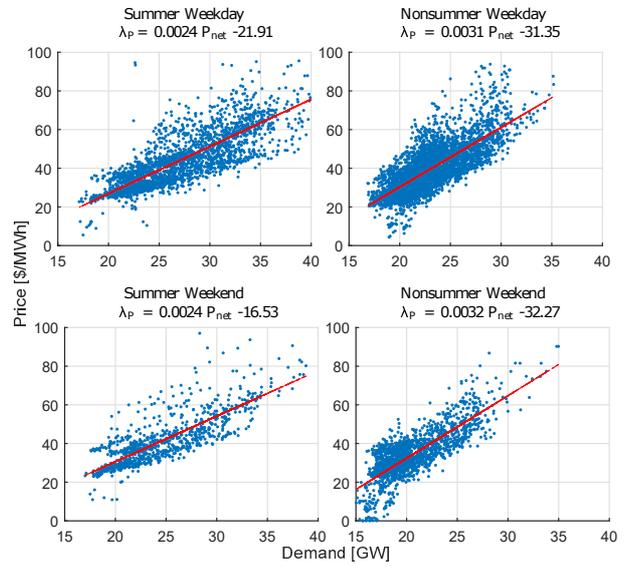}
\centering
\caption{CAISO price model fits as a function of net demand for each demand type (title).}
\label{figure:loadpricefitting}
\end {figure}

The MCP models developed for each demand cluster are shown in Fig. \ref{figure:loadpricefitting}. 

\begin{figure}[ht]
\includegraphics[width=.9\columnwidth]{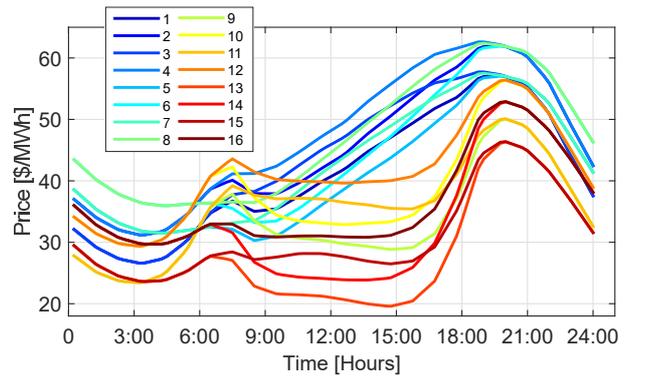}
\centering
\caption{Market clearing price (MCP) $\lambda_P$ for different scenarios (Table \ref{table:scenarios}) for the first year.
}\label{figure:Pricing}
\end {figure}

\begin{figure*}[ht]
\includegraphics[width=1.8\columnwidth,height=.38\columnwidth]{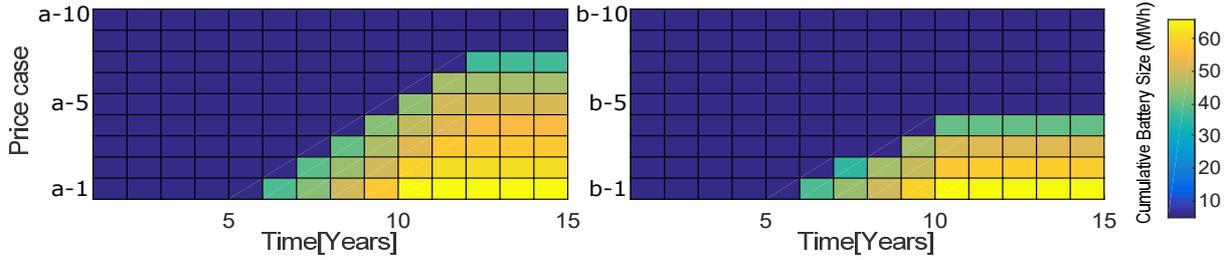}
\centering
\caption{BESS installation schedule by year for BESS price scenarios (a) on the left and (b) on the right.
}
\label{figure:Batplan}
\end {figure*}

The CAISO MCP $\lambda_P$ derived from Eq. \fref{price1} is the main input to the market model to determine the size and the daily operation of the BESS. The parameter $\lambda_P$ varies based on the different scenarios as shown in Fig.~\ref{figure:Pricing}. The highest market price is associated with low solar, low wind and summer weekday demand (scenario 4). Conversely, the scenario with non-summer weekend demand, high solar, and high wind results in the smallest $\lambda_P$. Negative pricing may appear also as a result of the assumed inability to curtail renewables; after renewable generators lose their protected ``must-take'' status, they will be curtailed in such a situation to avoid negative pricing. In this case, it is cheaper to temporarily pay market participants to take power rather than turning off base-load power plants. 

\subsection{Case Study}
The case study uses CAISO demand data and the location marginal pricing (LMP) node (UCM\_6\_N001) located at (32$^{\circ}$53'00.9"N 117$^{\circ}$13'21.2"W) which is the trading node containing UC San Diego. 
The simulated case has a peak market demand around 45~GW and a low demand (base-load) around 18~MW. The 2015 utility scale solar and wind peaks are around 5.7 and 2~GW respectively. 
Clear solar days are assumed to occur 30\% of the time, overcast days 70\%, and high wind 40\% of the time and low wind 60\%. Different demand scenarios follow the calendar with 96 and 36 days for summer (May 1 to Oct 31) weekday and weekend, respectively and 165 and 68 days for non-summer weekday and weekend, respectively. By that, the $Pr$ shown in \fref{Pr} is given by $Pr(Scen 1)=Pr(\left\{HS, HW, SWD\right\})=0.30\cdot 0.40\cdot 96/365=3.2\% $. The time step for all data and optimization schedules is 15 minutes.





\begin{figure}[ht]
\includegraphics[width=.85\columnwidth]{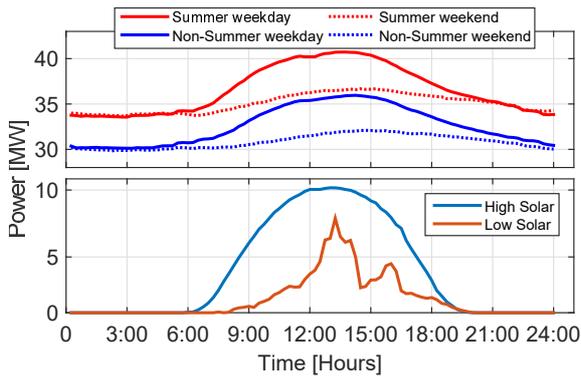}
\centering
\caption{Microgrid demand (top) and solar (bottom) profiles.} 
\label{figure:microdemandandsolar}
\end {figure} 

The characteristics of the microgrid conform to the UC San Diego microgrid. Similar to CAISO, demand is split into summer and non-summer weekday and weekend based on actual demand data collected from campus with demand peaks of (42, 36, 35, and 32)~MW for summer weekday, summer weekend, non-summer weekday, and non-summer weekend and base-load of 34~MW for both summer scenarios and 30~MW for non-summer as shown in Fig.~\ref{figure:microdemandandsolar} (top plot). Those profiles are matched with the existing market scenarios in Table~\ref{table:scenarios}. Microgrid generation is 20 MW from gas turbines and solar power of peak-to-peak ratio 
is 10~MW and both \textit{high} and \textit{low} solar clusters are shown in Fig.~\ref{figure:microdemandandsolar} (bottom plot). Noted here that the \textit{low} solar profile for the microgird is more intermittent compared to the market case because of the geospatial effect. The maximum allowable power demand from the grid ($P_L$) is 45~MW. Since microgrid energy sales to the market are not permitted, overgeneration would have to be curtailed. 


\begin{figure}[ht]
\includegraphics[width=.9\columnwidth]{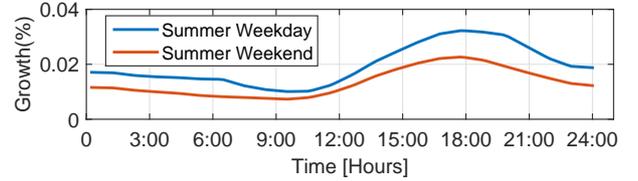}
\centering
\caption{The market annual demand growth profiles for summer.} 
\label{figure:marketdemandg}
\end {figure}

The growth rates of solar and demand are both $3\%$ for the microgrid. The growth of solar and wind are both assumed to be $7\%$ for the microgrid. To compute annual demand growth profile, CAISO demand data in 2013--2015 is used. The market annual demand growth profiles for non-summer scenarios are assumed scalar and equal to  2\% while those for summer scenarios are shown in Fig. \ref{figure:marketdemandg}. 


The BESS pricing cases ($J_{b,t}$ in Eq. \fref{min}), which all include government subsidies and incentives, are divided into two categories named $a$ and $b$.
Each category includes 10 cases (called $a\text{-}1,\cdots,a\text{-}10, b\text{-}1,\cdots,b\text{-}10$) which start from a price value between 117~\$/kWh and 175.5~\$/kWh. The price functions of all cases in category $a$ converge to 100 \$/kWh within 15 years \citep{nykvist2015rapidly} while those of the cases in category $b$ decay with a constant rate of 1\% every year. The life cycle of BESS is assumed 10 years for all cases.


\section{Numerical Results}\label{sec:results}

Fig.~\ref{figure:Batplan} shows optimization results for the BESS installation by year. On the $y$ axis prices increase from lower to higher. On the $x$ axis prices decrease from left to right as the years progress. For both cases no installation was applied before the year 5 (Y5) but as the BESS prices drop faster in case a compared to case b the installation went up to case $a\text{-}8$ compared to case $b\text{-}4$. 
The yearly installation plan of the BESS results to be large at one year followed by smaller installation few years before and after. After case $a\text{-}8$ and case $b\text{-}4$ no installations have resulted.

\begin{figure}[ht]
\includegraphics[width=.9\columnwidth]{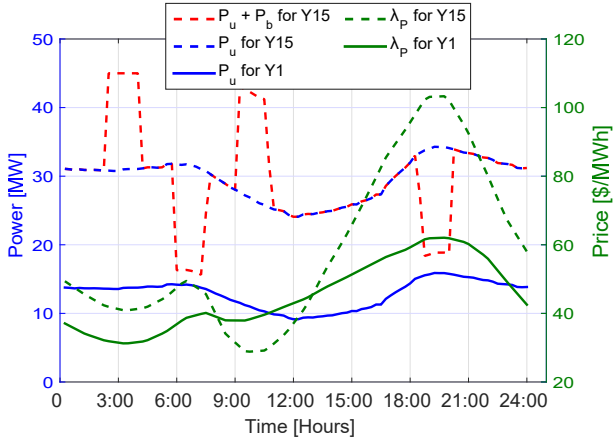}
\centering
\caption{First and last year microgrid demand and storage profile (left axis) and price (right axis). Data are for the first scenario \{SWD,HS,HW\}. During Y1 no BESS installation was present and no data is plotted. 
} 
\label{figure:Sampleday}
\end {figure}



Fig.~\ref{figure:Sampleday} shows a sample day of price case $a\text{-}8$ as $P_u$, $P_b$ and $\lambda_P$, where microgrid power purchase with storage is presented in the left axis and  $\lambda_P$ is presented in the right axis  for Y1 and Y15. Adding BESS in the microgrid changes the behavior of the microgrid demand with a new peak between 1 and 3~AM during a market price depression. Specifying $P_L = 45$~MW  limits the charging of the BESS from hour 1AM to 3AM to not exceed a total microgrid demand of 45~MW. Comparing $\lambda_P$ in Y1 and Y15 shows the effect of larger solar for a midday price minimum and stronger peak demand growth in the evening. Therefore the price curve shows two peaks and a pronounced evening peak in Y15. The price curve triggers two BESS cycles per day to leverage the margins between minima and peaks.

\section{Conclusions}
As the need for Battery Energy Storage System (BESS) is increasing to cope with intermittent energy production from renewable energy sources, optimal BESS sizing from an economical perspective is a challenging problem. 
To optimally size the BESS in the design stage of a microgrid, the trade-off between BESS cost, energy bill savings, and lifetime must be taken into account. 
Using a stochastic optimization approach we optimally size and schedule a BESS in a microgrid based on market energy pricing. Variability of wind and solar energy resources, the variability of energy demand, and a dynamic market price model that considers feedback from microgrid energy decisions are considered. Decreasing BESS costs over time are also modeled. The modeling framework contains significant flexibility and realism for microgrid planning. 

Assuming the market clearance price model and net demand forecasts in our case study for CAISO and a particular trading node, the results show that a microgrid can start saving money with wholesale market energy trading once BESS prices drop below \$150/kWh. The operational scheduling of BESS is not targeted to shave the microgrid peak but rather to profit from wholesale energy cost margins. Additional constraints could be added to achieve a hybrid between local and market objectives.

\bibliography{ifacconf}

\end{document}